\documentclass[smallcondensed, nospthms]{svjour3}
\journalname{Mathematical Programming, Ser. B}
\smartqed  
\usepackage{graphicx}
\usepackage[nohug]{diagrams}
\diagramstyle[labelstyle=\scriptstyle]
\usepackage{fix-cm}
\usepackage{latexsym}
\usepackage{amssymb}
\usepackage[center]{caption}
\usepackage{mathptmx}
\newcommand{\nc}{\newcommand}
\nc{\nt}{\newtheorem} 
\nt{thm}{Theorem}[section]
\nt{cor}[thm]{Corollary} 
\nt{prop}[thm]{Proposition}
\nt{lem}[thm]{Lemma} 
\nt{defn}[thm]{Definition}
\nt{exa}[thm]{Example} 
\nt{rem}[thm]{Remark}
\nt{ass}[thm]{Assumption}
\nt{alg}[thm]{Algorithm}
\nt{con}[thm]{Conjecture} 
\nc{\ip}[2]{\mbox{$\langle #1,#2 \rangle$}} 
\nc{\pf}{\noindent{\bf Proof\ \ }}
\nc{\finpf}{\hfill{$\Box$}\linespace}
\nc{\linespace}{\vspace{\baselineskip} \noindent} 
\nc{\R}{{\bf R}}
\nc{\cl}{\mbox{\rm cl}\,} 
\nc{\conv}{\mbox{\rm conv}\,} 
\nc{\gph}{\mbox{\rm gph}\,} 
\nc{\epi}{\mbox{\rm epi}\,} 
\nc{\dom}{\mbox{\rm dom}\,} 
\nc{\st}{\mbox{\rm s.t.}\,} 
\nc{\ri}{\mbox{\rm ri}\,} 
\nc{\inter}{\mbox{\rm int}\,} 
\nc{\para}{\mbox{\rm par}\,}

\newenvironment{myequation}{\begin{equation}}{\end{equation}}

\newenvironment{myeqnarray*}{\begin{eqnarray*}}{\end{eqnarray*}}
\nc{\bmye}{\begin{myequation}} \nc{\emye}{\end{myequation}}

\def\tto{\;{\lower 1pt \hbox{$\rightarrow$}}\kern -12pt
           \hbox{\raise 2.8pt \hbox{$\rightarrow$}}\;}

\usepackage{enumerate}
\usepackage{graphicx}

\usepackage{url}
\begin{document}

\title{Semi-algebraic functions have small subdifferentials}

\author{D. Drusvyatskiy         \and
        A.S. Lewis
}



\institute{D. Drusvyatskiy \at
    School of Operations Research and Information Engineering,
    Cornell University,
    Ithaca, New York, USA; 
    \email{dd379@cornell.edu}.   
    Work of Dmitriy Drusvyatskiy on this paper has been partially supported by the NDSEG grant from the Department of Defense.        
           \and
  A.S. Lewis \at
  School of Operations Research and Information Engineering,
  Cornell University,
  Ithaca, New York, USA;
  {\tt http://people.orie.cornell.edu/{\raise.17ex\hbox{$\scriptstyle\sim$}}aslewis/}.
  Research supported in part by National Science Foundation Grant DMS-0806057.
}

\date{Received: date / Accepted: date}
\maketitle

\begin{abstract}
We prove that the subdifferential of any semi-algebraic extended-real-valued function on $\R^n$ has $n$-dimensional graph.  We discuss consequences for generic semi-algebraic optimization problems.
\end{abstract}

\section{Introduction}
A principal goal of variational analysis is the search for generalized critical points of nonsmooth functions $f \colon \R^n \to \R$.  For example, given a locally Lipschitz function $f$, we might be interested in points $x \in \R^n$ having zero in the ``Clarke generalized gradient'' (or ``subdifferential'')  
$\partial_c f(x)$, a set consisting of convex combinations of limits of gradients of $f$ at points near 
$x$ \cite{clarke}.  

Adding a linear perturbation, we might seek critical points of the function
$x \mapsto f(x) - v^T x$ for a given vector $v \in \R^m$, or, phrased in terms of the graph of the subdifferential mapping $\partial_c f$, solutions to the inclusion
\[
(x,v) \in \mbox{gph}\,\partial_c f.
\]
More generally, given a smooth function $G \colon \R^m \to \R^n$, we might be interested in solutions 
$(x,y) \in \R^m \times \R^n$ to the system
\bmye \label{inclusion}
(G(x),y) \in \mbox{gph}\,\partial_c f ~~\mbox{and}~~  \nabla G(x)^* y = v
\emye
(where $*$ denotes the adjoint).  Such systems arise naturally when we seek critical points of the composite function 
$x \mapsto f(G(x)) - v^T x$.

Generalized critical points of {\em smooth} functions $f$ are, of course, simply the critical points in the classical sense.  However, the more general theory is particularly interesting to optimization specialists, because critical points of continuous convex functions are just minimizers \cite[Proposition 8.12]{VA}, and more generally, for a broader class of functions (for instance, those that are Clarke regular \cite{clarke}), a point is critical exactly when the directional derivative is nonnegative in every direction.

The system (\ref{inclusion}) could, in principle, be uninformative if the graph $\mbox{gph}\,\partial_c f$ is large.  In particular, if the dimension (appropriately defined) of the graph is larger than $n$, then we could not typically expect the system to be a very definitive tool, since it involves $m+n$ variables constrained by only $m$ linear equations and the inclusion. Such examples are not hard to construct:  indeed, there exists a function $f \colon \R \to \R$ with Lipschitz constant one and with the property that its Clarke subdifferential is the interval $[-1,1]$ at every point \cite{example}.  Alarmingly, in a precise mathematical sense, this property is actually typical for such functions \cite{gen_lip}.

Optimization theorists often consider subdifferentials that are smaller than Clarke's, the ``limiting'' subdifferential $\partial f$  being a popular choice \cite{VA,CLSW,Mord_1,Borwein-Zhu}.  However, the Clarke subdifferential can be easier to approximate numerically (see \cite{BLO}), and in any case the potential difficulty posed by functions with large subdifferential graphs persists with the limiting subdifferential \cite{badlim}.

Notwithstanding this pathology, concrete functions $f \colon \R^n \to \R$ encountered in practice have subdifferentials $\partial_c f$ whose graphs are, in some sense, small and this property can be useful, practically.  For instance, Robinson \cite{rob} considers algorithmic aspects of functions whose subdifferential graphs are everywhere locally Lipschitz homeomorphic to an open subset of $\R^n$.  As above, dimensional considerations suggest reassuringly that this property should help the definitive power of critical point systems like (\ref{inclusion}), and Robinson furthermore argues that it carries powerful computational promise. An example of the applicability of Robinson's techniques is provided by Minty's theorem, which states that the graph of the subdifferential of a proper, lower semicontinuous, convex function $f\colon\R^n\to\overline{\R}$ is Lipschitz homeomorphic to $\R^n$ \cite{minty}. 

When can we be confident that a function has a subdifferential graph that is, by some definition, small?  The study of classes of functions that are favorable for subdifferential analysis, in particular excluding the pathological examples above, is well-developed.  The usual starting point is a unification of smooth and convex analysis, arriving at such properties as amenability \cite[Chapter 10.F.]{VA}, prox-regularity \cite{prox_reg}, and cone-reducibility \cite[Section 3.4.4]{Bon_Shap}. Using Minty's theorem, Poliquin and Rockafellar \cite{prox_reg} showed that prox-regular functions, in particular, have small subdifferentials in the sense of Robinson.  Aiming precisely at a class of functions with small subdifferentials (in fact minimal in the class of upper semicontinuous mappings with nonempty compact convex images), \cite{Bor-Moors} considers ``essential strict differentiability''.

In this work we take a different, very concrete approach.  We focus on the dimension of the subdifferential graph, unlike the abstract minimality results of \cite{Bor-Moors}, but we consider the class of {\em semi-algebraic} functions---those functions whose graphs are semi-algebraic, meaning composed of finitely-many sets, each defined by finitely-many polynomial inequalities---and prove that such functions have small subdifferentials in the sense of dimension: the Clarke subdifferential has $n$-dimensional graph.  This result subsumes neither the simple case of a smooth function, nor the case of a convex function, neither of which is necessarily semi-algebraic. Nonetheless, it has a certain appeal:  semi-algebraic functions are common, they serve as an excellent model for ``concrete'' functions in variational analysis \cite{tame_opt}, and in marked contrast with many other classes of favorable functions, such as amenable functions, they may not even be Clarke regular. Furthermore, semi-algebraic functions are easy to recognize (as a consequence of the Tarski-Seidenberg theorem on preservation of semi-algebraicity under projection). For instance, observe that the spectral radius function on $n\times n$ matrices is neither Lipschitz nor convex, but it is easy to see that it is semi-algebraic. 

To illustrate our results, consider the critical points of the function $x \mapsto f(x) - v^T x$ for a semi-algebraic function $f \colon \R^n \to [-\infty,+\infty]$.  As a consequence of the subdifferential graph being small, we show that for a {\em generic} choice of the vector $v$, the number of critical points is finite.  More precisely, there exists a number $N$, and a semi-algebraic set $S \subset \R^n$ of dimension strictly less than $n$, such that for all vectors $v$ outside $S$, there exist at most $N$ critical points.
A result of a similar flavor can be found in \cite{crit_val}, where criticality of so called ``constraint systems'' is considered. Specifically, \cite{crit_val} shows that if a semi-algebraic constrained minimization problem is ``normal'', then it has only finitely many critical points.  Furthermore, it is shown that normality is a generic property. To contrast their approach to ours, we should note that \cite{crit_val} focuses on perturbations to the constraint structure, whereas we address linear perturbations to the function itself. 

To be concrete, we state our results for semi-algebraic functions.  Analogous results, with essentially identical proofs, hold for functions definable in an ``o-minimal structure'' and, more generally, for ``tame'' functions. (In the case of tame functions, ``finiteness'' of critical points should be replaced by ``local isolation'' in Proposition 4.3 and Corollaries 4.4,~\ref{cor:comp_unique2},~\ref{cor:gen_unique2}.) In particular, our results hold for globally subanalytic functions, discussed in \cite{Shiota}. For a quick introduction to these concepts in an optimization context, see \cite{tame_opt}.

\section{Preliminaries}
\subsection{Variational Analysis}
In this section, we summarize some of the fundamental tools used in variational analysis and nonsmooth optimization.
We refer the reader to the monographs of Rockafellar-Wets \cite{VA}, Borwein-Zhu \cite{Borwein-Zhu},
Mordukhovich \cite{Mord_1,Mord_2}, and Clarke-Ledyaev-Stern-Wolenski \cite{CLSW}, for more details.  Unless otherwise stated, we follow the terminology and notation of \cite{VA}.

Consider the extended real line $\overline{\R}:=\R\cup\{-\infty\}\cup\{+\infty\}$. We say that an extended-real-valued function is proper if it is never $\{-\infty\}$ and is not always $\{+\infty\}$.  

For a function $f\colon\R^n\rightarrow\overline{\R}$, we define the {\em domain} of $f$ to be $$\mbox{\rm dom}\, f:=\{x\in\R^n: f(x)<+\infty\},$$ and we define the {\em epigraph} of $f$ to be $$\mbox{\rm epi}\, f:= \{(x,r)\in\R^n\times\R: r\geq f(x)\}.$$

A {\em set-valued mapping} $F$ from $\R^n$ to $\R^m$, denoted by $F\colon\R^n\rightrightarrows\R^m$, is a mapping from $\R^n$ to the power set of $\R^m$. Thus for each  point $x\in\R^n$, $F(x)$ is a subset of $\R^m$. For a set-valued mapping $F\colon\R^n\rightrightarrows\R^m$, we define the {\em domain} of $F$ to be $$\mbox{\rm dom}\, F:=\{x\in\R^n:F(x)\neq\emptyset\},$$ and we define the {\em graph} of $F$ to be $$\mbox{\rm gph}\, F:=\{(x,y)\in\R^n\times\R^m:y\in F(x)\}.$$

\begin{defn}
{\rm Consider a set-valued mapping $F\colon\R^n\rightrightarrows\R^m$. 
\begin{enumerate}
\item $F$ is {\em outer semicontinuous} at a point ${\bar x}\in\R^n$ if for any sequence of points $x_r\in\R^n$ converging to $\bar{x}$ and any sequence of points $y_r\in F(x_r)$ converging to $\bar{y}$, we must have $\bar{y}\in F(\bar{x})$. 
\item $F$ is {\em inner semicontinuous} at $\bar{x}$ if for any sequence of points $x_r\in\R^n$ converging to $\bar{x}$ and any point $\bar{y}\in F(\bar{x})$, there exists a sequence $y_r\in\R^m$ converging to $\bar{y}$ such that $y_r\in F(x_r)$ for all $r$. 
\end{enumerate}
If both properties hold, then we say that $F$ is {\em continuous} at $\bar{x}$.
}
\end{defn}

\begin{defn}
{\rm Consider a set $S\subset\R^n$ and a point $\bar{x}\in S$. The {\em regular normal cone} to $S$ at $\bar x$, denoted 
$\hat N_S(\bar x)$, consists of all vectors $v \in \R^n$ such that $$\langle v,x-\bar{x} \rangle \leq o(|x-\bar{x}|) \textrm{ for }x\in S,$$ 
where we denote by $o(|x-\bar{x}|) \textrm{ for }x\in S$ a term with the property that $$\frac{o(|x-\bar{x}|)}{|x-\bar{x}|}\rightarrow 0$$ when $x\stackrel{S}{\rightarrow} \bar{x}$ with $x\neq\bar{x}$.
} 
\end{defn}

Given a closed set $S$, the mapping $x\mapsto\hat{N}_S(x)$ does not necessarily have a closed graph. To correct for that, the following definition is introduced.
\begin{defn}
{\rm Consider a set $S\subset\R^n$ and a point $\bar{x}\in S$.  The {\em limiting normal cone} to $S$ at $\bar{x}$, denoted $N_S(\bar{x})$, consists of all $v\in\R^n$ such that there are sequences $x_r\stackrel{S}{\rightarrow} \bar{x}$ and $v_r\rightarrow v$ with $v_r\in\hat{N}_S(x_r)$.}
\end{defn}

For a set $S\subset\R^n$, we denote its topological closure by $\mbox{\rm cl}\, S$ and its convex hull by $\mbox{\rm conv}\, S$. 

\begin{defn}
{\rm Consider a set $S\subset\R^n$ and a point $\bar{x}\in S$. The {\em Clarke normal cone} to $S$ at $\bar{x}$, denoted $N_S^c(\bar{x})$, is defined by $$N_S^c(\bar{x})=\textrm{cl conv }N_S(\bar{x}).$$}
\end{defn}

We summarize some simple facts about normal cones that we will need.
\begin{thm}\label{thm:norm_prop}
Consider a set $S\subset\R^n$ and a point $\bar{x}\in S$.
\begin{enumerate}
\item $\hat{N}_S(\bar{x})\subset N_S(\bar{x}) \subset N_S^c(\bar{x})$.
\item $N_S(\bar{x})$, $\hat{N}_S(\bar{x})$, and $N_S^c(\bar{x})$ are closed cones. $\hat{N}_S(\bar{x})$ and $N_S^c(\bar{x})$ are, in addition, convex. 
\item For a set $F\subset\R^n$ containing $\bar{x}$ such that $S\subset F$, we have $\hat{N}_F(\bar{x})\subset \hat{N}_S(\bar{x})$.
\end{enumerate}
\end{thm}

\begin{defn}[Clarke regularity of sets]
{\rm A set $S\subset\R^n$ is said to be {\em Clarke regular} at a point $\bar{x}\in S$ if it is locally closed at $\bar{x}$ and every limiting normal vector to $S$ at $\bar{x}$ is a regular normal vector, that is $N_S(\bar{x})=\hat{N}_S(\bar{x})$.}
\end{defn}

Given any set $S\subset\R^n$ and a mapping $f\colon S\to \widetilde{S}$, where $\widetilde{S}\subset\R^m$, we say that $f$ is {\em smooth} if for each point $\bar{x}\in S$, there is a neighborhood $U$ of $\bar{x}$ and a $\bf{C}^1$ mapping $\hat{f}\colon \R^n\to\R^m$ that agrees with $f$ on $S\cap U$. If a smooth function $f$ is bijective and its inverse is also smooth, then we say that $f$ is a {\em diffeomorphism}. 

What we call smooth is usually referred to as $\bf{C}^1$ smooth. Since in this work we will not need higher order of smoothness, no ambiguity should arise.

\begin{defn}[{\cite[Proposition 8.12]{Lee}}]
{\rm Consider a set $M\subset\R^n$. We say that $M$ is a {\em manifold} (or ``embedded submanifold'') of dimension $r$ if for each point $\bar{x}\in M$, there is an open neighborhood $U$ around ${\bar{x}}$ such that $M\cap U=F^{-1}(0)$, where $F\colon U\to\R^{n-r}$ is a smooth map with $\nabla F(\bar{x})$ of full rank. In this case, we call $F$ a {\em local defining function} for $M$ around $\bar{x}$.  
}
\end{defn}

\begin{thm}[{\cite[Example 6.8]{VA}}]\label{thm:clarke_man}
If $M$ is a manifold, then for every point $\bar{x}\in M$, the manifold $M$ is Clarke regular at $\bar{x}$ and $N_M(\bar{x})$ is equal to the normal space to $M$ at $\bar{x}$, in the sense of differential geometry.
\end{thm}


Normal cones allow us to study geometric objects. We now define subdifferentials, which allow us to analyze behavior of functions.
\begin{defn} 
{\rm Consider a function $f\colon\R^n\rightarrow\overline{\R}$ and a point $\bar{x}\in\R^n$ where $f$ is finite. The {\em regular}, {\em limiting}, and {\em Clarke subdifferentials} of $f$ at $\bar{x}$, respectively, are defined by 
$$\hat{\partial}f(\bar{x})= \{v\in\R^n: (v,-1)\in \hat{N}_{\mbox{{\scriptsize {\rm epi}}}\, f}(\bar{x},f(\bar{x}))\},$$  
$$\partial f(\bar{x})= \{v\in\R^n: (v,-1)\in N_{\mbox{{\scriptsize {\rm epi}}}\, f}(\bar{x},f(\bar{x}))\},$$  
$$\partial_c f(\bar{x})= \{v\in\R^n: (v,-1)\in N^c_{\mbox{{\scriptsize {\rm epi}}}\, f}(\bar{x},f(\bar{x}))\}.$$}  
\end{defn}

\noindent For $x$ such that $f(x)$ is not finite, we follow the convention that $\hat{\partial}f(x)=\partial f(x)=\partial_c f(x)=\emptyset$.

\begin{defn}[Subdifferential regularity]
{\rm A function $f\colon\R^n\to\overline{\R}$ is called {\em subdifferentially regular} at $\bar{x}$ if $f(\bar{x})$ is finite and $\mbox{\rm epi}\, f$ is Clarke regular at $(\bar{x},f(\bar{x}))$ as a subset of $\R^n\times\R$.}
\end{defn}

\begin{thm}[{\cite[Exercise 8.8, Corollary 10.9]{VA}}]\label{thm:sum} 
Consider the function $h=f+g$, where $f\colon\R^n\to\overline{\R}$ is finite at $\bar{x}$ and $g\colon\R^n\to\overline{\R}$ is smooth on a neighborhood of $\bar{x}$. Then we have $$\hat{\partial}h(\bar{x})= \hat{\partial}f(\bar{x})+\nabla g(\bar{x}),~~ \partial h(\bar{x})=\partial f(\bar{x})+\nabla g(\bar{x}).$$ Furthermore, $h$ is subdifferentially regular at $\bar{x}$ if and only if $f$ is subdifferentially regular at $\bar{x}$.
\end{thm}

For a set $S\subset\R^n$, we define $\delta_S\colon\R^n\to\overline{\R}$ to be a function that is $0$ on $S$ and $+\infty$ elsewhere. We call $\delta_S$ the {\em indicator function} of the set $S$. 

\begin{thm}[{\cite[Exercise 8.14]{VA}}]\label{thm:indic}
Consider the indicator function $\delta_S$ of a set $S\subset\R^n$. Then we have $$\partial\delta_S(\bar{x})=N_S(\bar{x}),~~ \hat{\partial}\delta_S(\bar{x})=\hat{N}_S(\bar{x}).$$ Furthermore, $\delta_S$ is subdifferentially regular at $\bar{x}$ if and only if $S$ is Clarke regular at $\bar{x}$.  
\end{thm}



\subsection{Semi-algebraic Geometry}
A {\em semi-algebraic} set $S\subset\R^n$ is a finite union of sets of the form $$\{x\in \R^n: P(x)=0, Q_1(x)<0,\ldots, Q_l(x)<0\},$$ where $P,Q_1,\ldots,Q_l$ are polynomials in $n$ variables. In other words, $S$ is a union of finitely many sets, each defined by finitely many polynomial equalities and inequalities. A map $F\colon\R^n\rightrightarrows\R^m$ is {\em semi-algebraic} if $\mbox{\rm gph}\, F\subset\R^{n+m}$ is a semi-algebraic set. Semi-algebraic sets enjoy many nice structural properties. We discuss some of these properties in this section. See the monographs of Basu-Pollack-Roy \cite{ARAG}, Lou van den Dries \cite{LVDB}, and Shiota \cite{Shiota}. For a quick survey, see the article of van den Dries-Miller \cite{DM} and the surveys of Coste \cite{Coste-semi,Coste-min}. Unless otherwise stated, we follow the notation of \cite{DM} and \cite{Coste-semi}. 

A fundamental fact about semi-algebraic sets is provided by the Tarski-Seidenberg Theorem \cite[Theorem 2.3]{Coste-semi}. It states that the image of any semi-algebraic set $S\subset\R^n$, under a projection to any linear subspace of $\R^n$, is a semi-algebraic set. From this result, it follows that a great many constructions preserve semi-algebraicity. In particular, for a semi-algebraic function $f\colon\R^n\to\overline{\R}$, it is easy to see that the set-valued mappings $\hat{\partial} f$, $\partial f$, and $\partial_c f$ are semi-algebraic. See for example \cite[Proposition 3.1]{tame_opt}. 

The most striking and useful fact about semi-algebraic sets is that they can be partitioned into finitely many semi-algebraic manifolds that fit together in a regular pattern. The particular stratification that we are interested in is defined below. 
\begin{defn}\label{defn:whit}
{\rm Consider a semi-algebraic set $Q$ in $\R^n$. A {\em Whitney stratification} of $Q$ is a finite partition of $Q$ into semi-algebraic manifolds $M_i$ (called strata) with the following properties: 
\begin{enumerate}
\item For distinct $i$ and $j$, if $M_i\cap \mbox{\rm cl}\,{M_j} \neq \emptyset$, then $M_i\subset\mbox{\rm cl}\, {M_j}\setminus M_j$.
\item For any sequence of points $(x_k)$ in a stratum $M_j$ converging to a point $x$ in a stratum $M_i$, if the corresponding normal vectors $y_k \in N_{M_j}(x_k)$ converge to a vector $y$, then $y \in N_{M_i}(x)$.
\end{enumerate} 
}
\end{defn}

Observe that property 1 of Definition~\ref{defn:whit} gives us topological information on how the strata fit together, while property 2 gives us control over how sharply the strata fit together. Property 1 is called the {\em frontier condition} and property 2 is called {\em Whitney condition (a)}. We should note that Whitney stratification, as defined above, is normally referred to as $\bf{C}^1$-Whitney stratification. Furthermore, Whitney condition (a) is usually stated somewhat differently. The equivalence is noted in \cite{IS}. One simple example of this type of a stratification to keep in mind throughout the discussion is the partition of a polytope into its open faces.

\begin{defn}
{\rm Given finite collections $\{B_i\}$ and $\{C_j\}$ of subsets of $\R^n$, we say that $\{B_i\}$ is {\em compatible} with $\{C_j\}$ if for all $B_i$ and $C_j$, either $B_i\cap C_j=\emptyset$ or $B_i\subset C_j$.}
\end{defn}

As discussed above, the following theorem is true.
\begin{thm}[{\cite[Theorem 4.8]{DM}}]\label{thm:whitney}
Let $Q, C_1,\ldots,C_l$ be semi-algebraic sets in $\R^n$. Then $Q$ admits a Whitney stratification that is compatible with $C_1,\ldots,C_l$.
\end{thm}

The notion of a stratification being compatible with some predefined sets might not look natural; in fact, it is crucial since this property enables us to construct refinements of stratifications


We will have occasion to use the following result.

\begin{thm}[{\cite[Theorem 4.8]{DM}}]\label{theorem:strat}
Consider a semi-algebraic set $S$ in $\R^n$ and a semi-algebraic map $f\colon S\rightarrow\R^m$. Let $\mathcal{A}$ be a finite collection of semi-algebraic subsets of $S$ and $\mathcal{B}$
a finite collection of semi-algebraic subsets of $\R^m$. Then there exists a Whitney stratification $\mathcal{A}'$ of $S$ that is compatible with $\mathcal{A}$ and a Whitney stratification $\mathcal{B}'$ of $\R^m$ compatible with $\mathcal{B}$ such that for every stratum $Q\in\mathcal{A}'$, we have that the restriction $f|_Q$ is smooth and $f(Q)\in\mathcal{B}'$.
\end{thm}

In particular, it follows that semi-algebraic maps are ``generically'' (in a sense about to be made clear) smooth.


\begin{defn}
{\rm Let $A\subset\R^n$ be a nonempty semi-algebraic set. Then we define the {\em dimension} of $A$, $\dim A$, to be the maximal dimension of a stratum in any Whitney stratification of $A$. We adopt the convention that $\dim \emptyset=-\infty$.}  
\end{defn}

It can be easily shown that the dimension does not depend on the particular stratification. See \cite[Chapter 4]{LVDB} for more details. 
\begin{thm}Let $A$ and $B$ be nonempty semi-algebraic sets in $\R^n$. Then the following hold.
\begin{enumerate}
\item If $A\subset B$, then $\dim A\leq \dim B$.
\item $\dim A=\dim \mbox{\rm cl}\,{A}$.
\item $\dim (\mbox{\rm cl}\,{A}\setminus A)< \dim A$.
\item If $f\colon A\rightarrow\R^n$ is a semi-algebraic mapping, then $\dim f(A)\leq \dim A$. If $f$ is one-to-one, then $\dim f(A)=\dim A$. In particular, semi-algebraic homeomorphisms preserve dimension.
\item $\dim A\cup B= \max\{\dim A, \dim B\}$.
\item $\dim A\times B=\dim A+\dim B$.
\end{enumerate}
\end{thm}

We will need the following simple proposition.
\begin{prop}\label{prop:local}
Consider a Whitney stratification $\{M_i\}$ of a semi-algebraic set $Q\subset\R^n$. Let $M_j$ be a stratum of maximal dimension. Then for any point $\bar{x}\in M_j$, there exists a neighborhood $B\subset\R^n$ around $\bar{x}$ so that $$B\cap Q= B\cap M_j.$$
\end{prop}
{\pf Assume otherwise. Then there is a sequence $x_r\in Q$ converging to $\bar{x}$ with $x_r\notin M_j$. Since there are finitely many strata, we can assume that the whole sequence is contained in some stratum $M$. It follows that $\bar{x}$ is a limit point of $M$. By the frontier condition of the Whitney stratification, it must be that $\dim M_j< \dim M$, which is a contradiction since the stratum $M_j$ was chosen to have maximal dimension.
}\qed

A set $U\subset\R^n$ is said to be ``generic'', if it is large in some precise mathematical sense, depending on context. Two popular choices are that of $U$ being a {\em full-measure} set, meaning its complement has Lebesgue measure zero, and that of $U$ being {\em topologically generic}, meaning it contains a countable intersection of dense open sets. In general, these notions are very different. However for semi-algebraic sets, the situation simplifies drastically. Indeed, if $U\subset\R^n$ is a semi-algebraic set, then the following are equivalent.
\begin{itemize}
\item [$\bullet$] $U$ is full-measure.
\item [$\bullet$] $U$ is topologically generic.
\item [$\bullet$] The dimension of $U^c$ is strictly smaller than $n$.  
\end{itemize} 

We will say that a certain property holds for a generic vector $v\in\R^n$ if the set of vectors for which this property holds is generic in the sense just described. Generic properties of semi-algebraic optimization problems will be discussed in Section~\ref{sec:conseq}.

\begin{defn}~\label{def:triv}
{\rm Let $A\subset\R^m$ be a semi-algebraic set. A continuous semi-algebraic mapping $p\colon A\rightarrow\R^n$ is {\em semi-algebraically trivial} over a semi-algebraic set $C\subset\R^n$ if there is a semi-algebraic set $F$ and a semi-algebraic homeomorphism $h\colon p^{-1}(C)\rightarrow C\times F$ such that $p|_{p^{-1}(C)}={\rm proj}\circ h$, or in other words the following diagram commutes:}
\begin{diagram}[height=1.7em]
p^{-1}(C) &\rTo^h   &C\times F\\
          &\rdTo_p  &\dTo_{\mbox{\scriptsize {\rm proj}}} \\
          &         &C
\end{diagram}
{\rm We call $h$ a {\em semi-algebraic trivialization} of $p$ over $C$.}
\end{defn} 

Henceforth, we use the symbol $\cong$ to indicate that two semi-algebraic sets are semi-algebraically homeomorphic. 

\begin{rem} \label{rmk:Hardt}
{\rm If $p$ is trivial over some semi-algebraic set $C$, then we can decompose $p|_{p^{-1}(C)}$ into a homeomorphism followed by a simple projection. Also, since the homeomorphism $h$ in the definition is surjective and $p|_{p^{-1}(C)}={\rm proj}\circ h$, it follows that $h(p^{-1}(c))= \{c\}\times F$ for any $c\in C$. Thus for any point $c\in C$, we have $p^{-1}(c)\cong F$ and $p^{-1}(C)\cong C\times p^{-1}(c)$.}
\end{rem}

The following is a simple example of semi-algebraic triviality.
\begin{exa}\label{exa:hardt}
{\rm We follow the notation of Definition~\ref{def:triv}. Consider the semi-algebraic function $p\colon \R\rightarrow\R$ defined by $p(x)=x^2$. Now consider the semi-algebraic mapping $$h\colon \R\setminus\{0\}\rightarrow \R_{++}\times \{\pm 1\}, \qquad x\mapsto(x^2,\mbox{sgn}\, x).$$
It is easy to check that $h$ is a semi-algebraic homeomorphism, and furthermore we have $p={\rm proj}\circ h$ when restricted to $\R\setminus\{0\}$. Thus $h$ is a semi-algebraic trivialization of $p$ over $\R_{++}$.   
}
\end{exa}

\begin{defn}
{\rm In the notation of Definition~\ref{def:triv}, a trivialization $h$ is {\em compatible} with a semi-algebraic set $B\subset A$ if there is a semi-algebraic set $H\subset F$ such that $h(B\cap p^{-1}(C))= C\times H$.}
\end{defn}

If $h$ is a trivialization over $C$ then, certainly, for any set $B\subset A$ we know $h$ restricts to a homeomorphism from $B\cap p^{-1}(C)$ to $h(B\cap p^{-1}(C))$. The content of the definition above is that if $p$ is compatible with $B$, then $h$ restricts to a homeomorphism between $B\cap p^{-1}(C)$ and $C\times H$ for some semi-algebraic set $H\subset F$. Here is a simple example.
\begin{exa}
{\rm 
Let the semi-algebraic functions $p$ and $h$ be as defined in Example~\ref{exa:hardt}. Now notice that $h(\R_{++}\cap p^{-1}(\R_{++}))= \R_{++}\times \{+1\}$. Thus $h$ is compatible with $\R_{++}$.
}
\end{exa}

The following result will be used extensively in the rest of this work. See \cite[Chapter 9, Theorem 1.2]{LVDB} for more details.
\begin{thm}[Hardt triviality]\label{theorem:Hardt}
Let $A\subset\R^n$ be a semi-algebraic set and $p\colon A\rightarrow\R^m$, a continuous semi-algebraic mapping. There is a finite partition of the image $p(A)$ into semi-algebraic sets $C_1,\ldots, C_k$ such that $p$ is semi-algebraically trivial over each $C_i$. Moreover, if $Q_1,\ldots,Q_l$ are semi-algebraic subsets of $A$, we can require each trivialization $h_i\colon p^{-1}(C_i)\rightarrow C_i\times F_i$ to be compatible with all $Q_j$.
\end{thm}

\begin{exa}
{\rm
Consider the following elaboration on Example~\ref{exa:hardt}. Let the semi-algebraic functions $p$ and $h$ be defined as in Example~\ref{exa:hardt}. We saw that $h$ is a semi-algebraic trivialization of $p$ over $\R_{++}$. Let $f\colon\{0\}\rightarrow \{0\}\times \{0\}$ be the zero map. Observe $f$ is a semi-algebraic trivialization of $p$ over $\{0\}$. Thus $\{\R_{++}, \{0\}\}$ is a partition of $p(\R)$ guaranteed to exist by Theorem~\ref{theorem:Hardt}. 

}
\end{exa}

Given a continuous semi-algebraic function $p$, Theorem~\ref{theorem:Hardt} states that we can partition the image of $p$ into semi-algebraic sets $C_1,\ldots,C_k$, so that for each index $i=1,\ldots,k$, the restricted mapping $p|_{p^{-1}(C_i)}$ has a very simple form. By applying Theorem~\ref{theorem:Hardt} to various naturally occurring mappings, many interesting results can be obtained. See \cite[Chapter 9]{LVDB} for more details. In particular, by applying this theorem to the projection map we can break up semi-algebraic sets into simple building blocks that have product structure and analyze each one separately. This type of reasoning leads to the following corollary.
\begin{cor}\label{cor:svh}
Let $F \colon \R^n \rightrightarrows \R^m$ be a semi-algebraic set-valued mapping. Then there exists a partition of the domain of $F$ into semi-algebraic sets $X_1,X_2,\ldots, X_k$ with the following properties:
\begin{enumerate}
\item For each index $i=1,2,\ldots k$, there exists a 
semi-algebraic set $Y_i \subset \R^m$ and a semi-algebraic homeomorphism 
$\theta_i \colon \mbox{\rm gph}\,F|_{X_i} \to X_i \times Y_i$ satisfying
\[
\theta_i(\{x\} \times F(x)) = \{x\} \times Y_i~~ \mbox{for all}~ x \in X_i.
\]
Consequently, for all $x \in X_i$, we have $F(x) \cong  Y_i$ and 
$$\mbox{\rm gph}\,F|_{X_i}\cong X_i\times F(x).$$
\item If in addition, $\widetilde{F} \colon \R^n \rightrightarrows \R^m$ is another semi-algebraic set-valued mapping with $\widetilde{F}(x)\subset F(x)$, then we may also require that for each index $i=1,2,\ldots,k$, there exists a semi-algebraic set $\widetilde{Y}_i \subset Y_i$, such that $\theta_i(\mbox{\rm gph}\,\widetilde{F}|_{X_i})=X_i \times \widetilde{Y}_i$. Consequently, for all $x \in X_i$, we have $\widetilde{F}(x) \cong  \widetilde{Y}_i$ and 
$$\mbox{\rm gph}\,\widetilde{F}|_{X_i}\cong X_i\times \widetilde{F}(x).$$
\end{enumerate}
\end{cor}
{\pf Assume that we are given semi-algebraic set-valued maps $F$ and $\widetilde{F}$ such that $\widetilde{F}(x)\subset F(x)$ for all $x\in\R^n$. If $\widetilde{F}$ was not given, proceed with the proof with $\widetilde{F}(x)=\emptyset$ for all $x\in\R^n$. Consider $\mbox{\rm gph}\,F\subset\R^n\times\R^m$. Let $p\colon \mbox{\rm gph}\,F\rightarrow\R^n$ be the projection onto the first $n$ coordinates. By applying Theorem~\ref{theorem:Hardt} to $p$, we get a partition of the domain of $F$ into semi-algebraic sets $X_1,X_2,\ldots, X_k$ such that $p$ is semi-algebraically trivial over each $X_i$ and each trivialization is compatible with $\mbox{\rm gph}\,\widetilde{F}$. Thus there exist semi-algebraic sets $Y_1,Y_2,\ldots,Y_k \subset \R^m$ and $\widetilde{Y_1},\widetilde{Y_2},\ldots,\widetilde{Y_k} \subset \R^m$ with $\widetilde{Y_i}\subset Y_i$, such that for each $i$, there is a semi-algebraic homeomorphism $\theta_i\colon p^{-1}(X_i)\rightarrow X_i\times Y_i$, where we have 
\begin{eqnarray}
\theta_i(\mbox{\rm gph}\,\widetilde{F}\cap p^{-1}(X_i))&=& X_i\times \widetilde{Y_i},\nonumber\\ 
\label{cor:svh:triv}{\rm proj}_{X_i}\circ \theta_i&=&p|_{p^{-1}(X_i)}.  
\end{eqnarray}
Observe that $p^{-1}(X_i)=\mbox{\rm gph}\,F|_{X_i}$ and since $\mbox{\rm gph}\, \widetilde{F}$ is contained in $\mbox{\rm gph}\, F$, it follows that $\mbox{\rm gph}\,\widetilde{F}\cap p^{-1}(X_i)=\mbox{\rm gph}\,\widetilde{F}|_{X_i}$. Thus to summarize, we have 
\begin{eqnarray}
\label{cor:svh:cong1}\mbox{\rm gph}\,F|_{X_i}&\cong & X_i\times Y_i,\\
\nonumber \mbox{\rm gph}\,\widetilde{F}|_{X_i}&\cong & X_i\times\widetilde{Y_i}. 
\end{eqnarray}
Finally, from (\ref{cor:svh:triv}) and (\ref{cor:svh:cong1}), it follows that for all points $x \in X_i$, we have $$\theta_i(\{x\}\times F(x))=\{x\}\times Y_i,$$ completing the proof.
}\qed 

The following proposition appears in \cite{BZ2,BZ}; as observed there, this result is an easy and important consequence of Theorem \ref{theorem:Hardt}, and even though we will not have occasion to use it in this work, we include it and its proof below as an elegant illustration. 

\begin{prop}\label{prop:gen_cont}
Let $F\colon\R^n\rightrightarrows\R^m$ be a semi-algebraic set-valued mapping. Then there exists a finite partition of the domain of $F$ into semi-algebraic sets $X_1,\ldots,X_k$, such that for each index $i=1,\ldots,k$, the restricted mapping $F|_{X_i}$ is inner semicontinuous. If in addition, the mapping $F$ is compact-valued, then we can also require the restricted mapping $F|_{X_i}$ to be outer semicontinuous for each index $i=1,\ldots,k$. (In fact, the partition guaranteed to exist by Corollary~\ref{cor:svh} is one such partition.)
\end{prop}
{\pf
Applying Corollary~\ref{cor:svh} to the mapping $F$, we get a finite partition of the domain of $F$ into semi-algebraic sets $X_1,\ldots,X_k$, so that, in particular, property 1 of the corollary holds. To see the inner semicontinuity of the restricted map $F|_{X_i}$, consider any point $(\bar x, \bar{y}) \in \mbox{gph}\,F|_{X_i}$, and any sequence of points $x_r \to \bar x$ in the set $X_i$.  We want to construct a sequence of points $y_r \in F(x_r)$ converging to $\bar{y}$.  Notice that $\theta_i(\bar x,\bar y) = (\bar{x}, \hat y)$ for some point $\hat y \in Y_i$.  Since 
$(x_r,\hat y) \to (\bar x,\hat y)$, we deduce 
$\theta_i^{-1}(x_r,\hat y) \to \theta_i^{-1}(\bar x,\hat y) = (\bar{x},\bar y)$.  But for each index $r$, we know $\theta_i^{-1}(x_r,\hat y) = (x_r,y_r)$ for some point $y_r \in F(x_r)$, so the result follows. 

Assume now that $F$ is compact-valued. Consider any point $\bar x\in X_i$ and any sequence of points $(x_r,y_r)\to(\bar x,\bar y)$, where $\bar y$ is some point in $\R^m$ and $y_r\in F(x_r)$ for each $r$. We want to argue that $\bar{y}$ is in $F(\bar x)$. Consider the sequence $(\bar x,{\rm proj}_{Y_i}(\theta_i(x_r,y_r)))$. Observe that this sequence is contained in $\{\bar{x}\}\times Y_i$, which is a compact set since it is homeomorphic to $F(\bar x)$. Thus, without loss of generality, we can assume that $(\bar x,{\rm proj}_{Y_i}(\theta_i(x_r,y_r)))$ converges to $(\bar x, \hat{y})$ for some point $\hat{y}\in Y_i$. So we have $$(x_r,y_r)=\theta_i^{-1}(x_r,{\rm proj}_{Y_i}(\theta_i(x_r,y_r)))\to\theta_i^{-1}(\bar x, \hat{y})\in\{\bar x\}\times F(\bar x).$$ By the uniqueness of the limit, we must have $\bar y\in F(\bar x)$.}\qed

As a consequence of Proposition~\ref{prop:gen_cont}, it follows that any semi-algebraic set-valued mapping $F\colon\R^n\to\R^m$ is generically inner semicontinuous. If, in addition, $F$ is compact-valued, then $F$ is generically continuous. In fact, we can do better. If we require the mapping $F$ just to be closed-valued, then we can still partition its domain into semi-algebraic sets $X_1,\ldots,X_k$, such that for each index $i=1,\ldots,k$, the restricted mapping $F|_{X_i}$ is continuous. To see this, we need the following theorem that appears in \cite[Theorem 5.55]{VA}, and is attributed to \cite{GC,Kur,SCS}. Recall that given a topological space $X$, a subset $A$ of $X$ is meager if it is a union of countably many nowhere dense subsets of X.
\begin{thm}[Kuratowski]\label{thm:kur}
Consider a set $X\subset\R^n$ and a closed-valued set-valued mapping $F\colon X\rightrightarrows\R^m$. Assume that $F$ is either outer semicontinuous or inner semicontinuous relative to $X$. Then the set of points where $F$ fails to be continuous relative to $X$ is meager in $X$. 
\end{thm}
It is easy to see that if a semi-algebraic set $S$ is meager in another semi-algebraic set $X$, then the dimension of $S$ is strictly less than the dimension of $X$ (see \cite{gen} for more details). 
\begin{prop}
Let $F\colon\R^n\rightrightarrows\R^m$ be a semi-algebraic closed-valued set-valued mapping. Then there exists a finite partition of the domain of $F$ into semi-algebraic sets $X_1,\ldots,X_k$, such that for each index $i=1,\ldots,k$, the restricted mapping $F|_{X_i}$ is continuous. 
\end{prop}
{\pf
Applying Proposition~\ref{prop:gen_cont} to the mapping $F$, we get a partition of the domain of $F$ into semi-algebraic sets $X_1,\ldots,X_k$, so that the restricted map $F|_{X_i}$ is inner semicontinuous. Fix some set $X_i$. Let $S^0:=X_i$ and let $S^1\subset X_i$ be the set of points at which $F|_{S^0}$ fails to be continuous. By Theorem~\ref{thm:kur}, it follows that $\dim S^1<\dim S^0$. Now by applying this argument inductively, we can create a sequence of semi-algebraic sets $S^0\supset\ldots\supset S^k$, for some integer $k$, such that the collection $\{S_j\setminus S_{j+1}\}^{k-1}_{j=0}$ is a partition of $X_i$ and $F$ is continuous when restricted to each $S_j\setminus S_{j+1}$. By applying this argument to all the sets $X_i$, for $i=1,\ldots,k$, we get the result.
}\qed
\begin{rem}
{\rm In fact, it is shown in Daniilidis-Pang \cite{Pang} that closed-valued semi-algebraic maps are generically strictly continuous (see \cite{VA} for the definition). Their proof of this rather stronger result requires more sophisticated tools.}
\end{rem}
Finally, we have the following result:
\begin{thm}[{\cite[Theorem 4.4]{DM}}]\label{theorem:conn_fin}
Let $A$ be a semi-algebraic subset of $\R^n\times\R^m$. There is an integer $\beta$ such that for every point $x\in\R^n$, the number of connected components of the set $A_x=\{y\in \R^m:(x,y)\in A\}$ is no greater than $\beta$. 
\end{thm}

The following is a simple special case of Theorem~\ref{theorem:conn_fin}. We record it here for convenience.
\begin{rem}\label{cor:fiber}
{\rm Let $F\colon\R^n\rightrightarrows\R^m$ be a semi-algebraic mapping. Applying Theorem~\ref{theorem:conn_fin} to $\gph F\subset\R^n\times\R^n$, we deduce that there is an integer $\beta$ such that for every $x\in\R^n$, the number of connected components of $F(x)$ is no greater than $\beta$.}  
\end{rem}

\section{Main Results}\label{sec:main}
\begin{defn}
{\rm
Consider a Whitney stratification $\mathcal{A}$ of a semi-algebraic set $Q\subset\R^n$. We define the {\em normal bundle} $N_{\mathcal{A}}$ associated with the stratification $\mathcal{A}$ to be the union of the normal bundles of each stratum, that is 
$$N_{\mathcal{ A}}=\bigcup_{M\in\mathcal{A}} \gph N_M =\bigcup_{M\in\mathcal{A}} \{(x,y)\in \R^n\times\R^n: x\in M,\, y\in N_M(x)\}.$$}
\end{defn}

In the definition above, since there are finitely many strata and for each stratum $M\in\mathcal{A}$, the semi-algebraic set $\gph N_M$ is $n$-dimensional, we deduce that the normal bundle $N_{\mathcal{A}}$ is a semi-algebraic set of dimension $n$.

\begin{prop}\label{prop:strat}
Consider a semi-algebraic set $Q \subset \R^n$ and suppose it admits a Whitney stratification $\mathcal{A}=\{M_i\}$.  Then for any stratum $M_i$ and any point $\bar{x}\in M_i$, the Clarke normal cone, $N^c_Q(\bar{x})$, is contained in the normal space, $N_{M_i}(\bar{x})$. Consequently, the inclusion $\gph N_Q^c\subset N_{\mathcal{A}}$ holds and so the graph of the Clarke normal cone has dimension no greater than $n$. 
\end{prop}
{\pf
Observe that for any stratum $M_j$, we have the inclusion $M_j\subset Q$. Hence for any point $x\in M_j$, the inclusion 
\begin{equation}\label{eq:reg_inc}
\hat{N}_Q(x)\subset\hat{N}_{M_j}(x)=N_{M_j}(x) 
\end{equation}
holds. Now fix some stratum $M_i$ and a point $\bar{x}\in M_i$. We claim that the limiting normal cone $N_Q(\bar{x})$ is contained in $N_{M_i}(\bar{x})$. To see this, consider a vector $v\in N_Q(\bar{x})$. By definition of the limiting normal cone, there exist sequences $(x_r)$ and $(v_r)$ such that $x_r\stackrel{Q}{\rightarrow} \bar{x}$ and $v_r\rightarrow v$ with $v_r\in\hat{N}_Q(x_r)$. Since there are finitely many strata, we can assume that there is some stratum $M_j$ such that the entire sequence $(x_r)$ is contained in $M_j$. 
From (\ref{eq:reg_inc}), we deduce $\hat{N}_Q(x_r)\subset N_{M_j}(x_r)$, and hence $v_r\in N_{M_j}(x_r)$. Therefore by Whitney condition (a), we have $v\in N_{M_i}(\bar{x})$. Since $v$ was arbitrarily chosen from $N_Q(\bar{x})$, we deduce $N_Q(\bar{x})\subset N_{M_i}(\bar{x})$ and thus $N_Q^{c}(\bar{x})=\textrm{cl conv }N_Q(\bar{x})\subset N_{M_i}(\bar{x})$, as we needed to show.} \qed 

Shortly, we will generalize this result to the graph of the Clarke subdifferential. We need the following simple result. We provide a proof for completeness.

\begin{prop}\label{prop:const_gen}
Let $A\subset\R^n$ be a semi-algebraic set and $p\colon A\to\R^m$, a continuous semi-algebraic mapping. Let $D$ be the image set of the mapping $p$. Then we have the inequality, 
$$\dim D+ \min_{x\in D}\dim p^{-1}(x) \leq\dim A\leq \dim D+\max_{x\in D}\dim p^{-1}(x).$$
In particular, if there exists an integer $k$ such that the set $p^{-1}(x)$ is $k$-dimensional for every point $x\in D$, then the equality, $$\dim A= \dim D +k,$$ holds.   
\end{prop}
{\pf Applying Theorem~\ref{theorem:Hardt} to the mapping $p$, we obtain a finite partition of the set $D$ into semi-algebraic sets $\{C_i\}$ such that  
$p^{-1}(C_i)\cong C_i\times p^{-1}(c)$, for any $c\in C_i$. 
Let $C_i$ be a partitioning set satisfying $\dim A = \dim p^{-1}(C_i)$, and let $c$ be any point in $C_i$. Then we have $$\dim A = \dim p^{-1}(C_i)=\dim C_i+\dim p^{-1}(c) \leq\dim D+\max_{x\in D}\dim p^{-1}(x).$$ Let $C_j$ be a partitioning set satisfying $\dim D = \dim C_j$, and let $c$ be any point in $C_j$. Then we obtain $$\dim A \geq \dim p^{-1}(C_j)=\dim C_j+\dim p^{-1}(c) \geq\dim D+\min_{x\in D}\dim p^{-1}(x),$$ as we needed to show.
}\qed

We record the following simple and intuitive corollary for reference.
\begin{cor}\label{prop:const}
Let $F\colon\R^n\rightrightarrows\R^m$ be a semi-algebraic set-valued mapping and let $D:=\dom F$. Then we have the inequality, 
$$\dim D+\min_{x\in D}\dim F(x)\leq\dim \gph F\leq \dim D+\max_{x\in D}\dim F(x).$$ In particular, if there exists an integer $k$ such that the set $F(x)$ is $k$-dimensional for every point $x\in D$, then the equality, $$\dim \gph F= \dim D +k,$$ holds.   
\end{cor}
{\pf This is an easy consequence of Corollary~\ref{cor:svh}. 
}\qed


\begin{thm}\label{thm:main}
Let $f\colon\R^n\rightarrow\overline{\R}$ be a semi-algebraic function. Then the graph of the Clarke subdifferential, $\mbox{\rm gph}\, \partial_c f$, has dimension no greater than $n$. 
\end{thm}
{\pf
Let $F:=\epi f$ and $$A:= \{(x,r,y)\in \R^{n}\times{\R}\times\R^{n+1}:((x,r),y)\in\gph N_F^c,\, r=f(x),\,  y_{n+1}< 0\}.$$ Using Proposition~\ref{prop:strat}, we see 
\begin{equation}\label{eq:prev_eq}
\dim A\leq \dim\gph N_F^c\leq n+1. 
\end{equation}
Consider the continuous semi-algebraic map 
$$\phi\colon A\rightarrow\R^n\times\R^n$$ $$(x,f(x),y)\mapsto (x,\pi\Big(\frac{y}{|y_{n+1}|}\Big)),$$ where $\pi\colon\R^{n+1}\to\R^n$ is the canonical projection onto the first $n$ coordinates. Observe that the image of $\phi$ is exactly the graph of the Clarke subdifferential $\partial_c f$. Furthermore, for any pair $(x,v)\in \gph \partial_c f$, we have $$\phi^{-1}(x,v)=\{x\}\times \{f(x)\}\times \R_{+}(v,-1),$$ and hence $\dim \phi^{-1}(c)=1$ for any point $c$ in the image of $\phi$. By Proposition~\ref{prop:const_gen}, we deduce $$\dim \gph \partial_c f+1=\dim A\leq n+1,$$ where the last inequality follows from (\ref{eq:prev_eq}). Hence, we obtain $\dim \gph \partial_c f\leq n$, as we needed to show.
}\qed

Shortly we will show that for a proper semi-algebraic function $f\colon\R^n\rightarrow\overline{\R}$, both $\mbox{\rm gph}\, \partial_c f$ and $\mbox{\rm gph}\, \partial f$ have dimension exactly equal to $n$. In the case that the domain of $f$ is full-dimensional, this fact is easy to show. The argument is as follows. By Theorem~\ref{theorem:strat}, the domain of $f$ can be partitioned into semi-algebraic manifolds $\{X_i\}$ such that $f|_{X_i}$ is smooth. Let $X_i$ be the manifold of maximal dimension. Observe that for $x\in X_i$, we have $\partial f(x)=\{\nabla f(x)\}$ and it easily follows that $\dim \mbox{\rm gph}\, \partial f|_{X_i}=n$. Thus we have $$n\leq\dim \mbox{\rm gph}\, \partial f\leq \dim\mbox{\rm gph}\, \partial_c f\leq n,$$ where the last inequality follows from Theorem~\ref{thm:main}, and hence there is equality throughout. The argument just presented no longer works when the domain of $f$ is not full-dimensional. A slightly more involved argument is required. We record the following simple observation for reference.

\begin{prop}\label{prop:boundary}
Consider a smooth manifold $M\subset \R^n$ and a smooth real-valued function
$f\colon M\to\R$. Define a function $h:\R^n\to\overline{\R}$ agreeing with $f$ on $M$ and equaling plus infinity elsewhere. Then  $h$  is subdifferentially regular throughout  $M$.  Furthermore, at any point $\bar{x}\in M$, we have $$\partial h(\bar{x})=N_M(\bar{x})+\nabla g(\bar{x}),$$ where $g\colon\R^n\to\R$ is any smooth function agreeing with $f$ on $M$ on a neighborhood of $\bar{x}$. Consequently, $\partial h(\bar{x})$ is nonempty with dimension $n-\dim M$.
\end{prop}
{\pf 
Observe that near the point $\bar{x}$, we have $h=\delta_M+g$. Combining Theorem~\ref{thm:clarke_man} and Theorem~\ref{thm:indic}, we have that the function $\delta_M$ is subdifferentially regular at $\bar{x}$. By Theorem~\ref{thm:sum}, it follows that $h$ is subdifferentially regular at $\bar{x}$ and $$\partial h(\bar{x})=\partial\delta_M(\bar{x})+\nabla g(\bar{x})=N_M(\bar{x})+\nabla g(\bar{x}),$$ as we needed to show. 
}\qed

\begin{thm}\label{thm:equation}
Let $f\colon\R^n\rightarrow\overline{\R}$ be a proper semi-algebraic function. Then the graphs of the regular, limiting, and Clarke subdifferentials have dimension exactly $n$.
\end{thm}
{\pf We know $$\dim \mbox{\rm gph}\, \hat{\partial} f\leq\dim \mbox{\rm gph}\, \partial f\leq\dim \mbox{\rm gph}\, \partial_c f\leq n,$$ where the last inequality follows from Theorem~\ref{thm:main}. Thus if we show that the dimension of $\mbox{\rm gph}\, \hat{\partial} f$ is no less than $n$, we will be done. With that aim, applying Theorem~\ref{theorem:strat} to the function $f$, we obtain a Whitney stratification $\{M_i\}$ of the domain of $f$ such that for every stratum $M_i$, the restriction $f|_{M_i}$ is smooth. Let $M_j$ be a stratum of $\mbox{\rm dom}\, f$ of maximal dimension.  

Now consider the function $h\colon\R^n\to\overline{\R}$, which agrees with $f$ on $M_j$ and is plus infinity elsewhere. By Proposition~\ref{prop:local}, the functions $h$ and $f$ coincide on a neighborhood of $\bar{x}$. Applying Proposition~\ref{prop:boundary}, we deduce that $f$ is subdifferentially regular at $\bar{x}$ and $\partial f(\bar{x})$ is nonempty with dimension $n-\dim M_j$. Since the point $\bar{x}$ was arbitrarily chosen from $M_j$, we deduce $\dim\hat{\partial} f(x)=n-\dim M_j$ for any point $x\in M_j$. Thus applying Corollary~\ref{prop:const} to the semi-algebraic set-valued map $\hat{\partial} f|_{M_j}$, we deduce $\dim \gph \hat{\partial} f|_{M_j} = \dim M_j+n-\dim M_j=n$, and hence the result follows.  
}\qed 

More refined, local versions of Theorem~\ref{thm:equation} are investigated in \cite{loc}. 

Theorem~\ref{thm:main} shows that for a semi-algebraic function $f\colon\R^n\to\overline{\R}$, the Clarke subdifferential $\partial_c f$ is small in a dimensional sense. If $f$ is also Lipschitz, it is small in another sense, that we now discuss: we relate our results to the notion of a minimal {\bf cusco} ({\bf c}onvex {\bf u}pper {\bf s}emicontinuous nonempty {\bf co}mpact valued set-valued mapping), introduced in \cite{Bor-Moors}. To that effect, consider a set $A\subset\R^n$ and a set-valued mapping $F\colon A\rightrightarrows\R^m$. We say that $F$ is {\em upper semicontinuous} at some point $\bar{x}\in A$ if every open set $U$ containing $F(\bar{x})$ also contains $F(z)$ for all points $z\in A$ close to $\bar{x}$. If a map is closed-valued and upper semicontinuous, then it is outer semicontinuous. On the other hand, if a map is outer semi-continuous and locally bounded, then it is upper semicontinuous. See \cite[Section 5.1.4]{Borwein-Zhu} for more details. In particular, for a Lipschitz function $f$, the Clarke subdifferential $\partial_c f$ is upper semi-continuous \cite[Proposition 2.1.5]{clarke}. The mapping $F\colon A\rightrightarrows\R^m$ is said to be a {\em cusco} if it is upper-semicontinuous on $A$ and $F(x)$ is a nonempty compact convex set for each point $x\in A$. A {\em minimal cusco} is a cusco, whose graph does not strictly contain the graph of any other cusco. Let $U\subset\R^n$ be an open set and consider a semi-algebraic locally Lipschitz function $f\colon U\to\R$. It follows by a direct application of \cite[Corollary 2.2]{Bor-Moors} and generic smoothness of $f$ that the the set-valued mapping $\partial_c f$ is, in fact, a minimal cusco.

It is tempting to think that in the semi-algebraic setting, the graph of an arbitrary minimal cusco should have small dimension. However, it is not hard to see that this is not the case. For instance, we will now exhibit a semi-algebraic minimal cusco $F\colon \R^3\rightrightarrows\R^3$, whose graph is $4$-dimensional. Thus semi-algebraic minimal cuscos with low dimensional graphs, such as the Clarke subdifferential of a semi-algebraic locally Lipschitz function $f$ defined on an open set, are somewhat special. 

To simplify notation, we let $[y<0,z<0]$ be an alias for the set $\{(x,y,z)\in \R^3:y<0,z<0\}$ and we reserve analogous notation for relaters `$>$' and `$=$'. 
Consider the semi-algebraic set-valued mapping $F\colon\R^3\rightrightarrows\R^3$, defined as follows 
$$F|_{[y>0, z>0]}=\{(0,0,0)\}, ~F|_{[y<0, z>0]}=\{(0,0,1)\},$$ $$F|_{[y<0,z<0]}=\{(0,1,0)\}, ~F|_{[y>0, z<0]}=\{(1,0,0)\},$$ 
$$F|_{[y>0, z=0]}=\mbox{\rm conv}\,\{(0,0,0),(1,0,0)\}, ~F|_{[y=0, z>0]}=\mbox{\rm conv}\,\{(0,0,1),(0,0,0)\},$$ $$F|_{[y<0,z=0]}=\mbox{\rm conv}\,\{(0,1,0),(0,0,1)\},~F|_{[y=0, z<0]}=\mbox{\rm conv}\,\{(1,0,0),(0,1,0)\},$$ $$F|_{[y=0,z=0]}= \mbox{\rm conv}\,\{(0,0,0),(0,0,1),(0,1,0),(1,0,0)\}.$$ It is easy to verify that $F$ is indeed a minimal cusco with a $4$-dimensional graph. In particular, Theorem~\ref{thm:equation} implies that $F$ is not the Clarke subdifferential mapping $\partial_c f$ for any semi-algebraic function $f\colon\R^3\rightrightarrows\overline{\R}$.    


\section{Consequences}\label{sec:conseq}
\begin{defn}
{\rm Consider a function $f\colon\R^n\rightarrow\overline{\R}$. We say that a point $x\in\R^n$ is {\em Clarke-critical} for the function $f$ if $0\in\partial_{c}f(x)$, and we call such a critical point $x$ {\em nondegenerate} if the stronger property $0 \in \mbox{ri}\, \partial_c f(x)$ holds.}
\end{defn}

Recall that for a proper convex function $f\colon\R^n\to\overline{\R}$ and a point $\bar{x}\in\mbox{\rm dom}\, f$, the subdifferentials $\hat{\partial}f(\bar{x})$, $\partial f(\bar{x})$, and $\partial_c f(\bar{x})$ all coincide and are equal to the convex subdifferential of $f$ at $\bar{x}$. So in this case, the notions of Clarke-criticality and Clarke-nondegeneracy reduce to more familiar notions from Convex Analysis. The importance of nondegeneracy for the sensitivity analysis of convex functions is well known:  in \cite{LOS}, for example, it is an underlying assumption for a pioneering conceptual approach to superlinearly convergent convex minimization algorithms. Consider the following largely classical theorem (see \cite[Proposition 1]{gen} and \cite{prep}).

\begin{thm}\label{thm:gen_conv}
Let $f\colon\R^n\to\overline{\R}$ be a proper convex function. Consider the collection of perturbed functions $h_v(x)=f(x)-\langle v,x\rangle$, parametrized by vectors $v\in\R^n$. Then for a full measure set of vectors $v\in\R^n$, the function $h_v$ has at most one minimizer, which furthermore is nondegenerate.
\end{thm}

Shortly, we will prove that a natural analogue of Theorem~\ref{thm:gen_conv} holds for arbitrary semi-algebraic functions, with no assumption of convexity. We will then reference an example of a locally Lipschitz function that is not semi-algebraic, and for which the conclusion of our analogous result fails, thus showing that the assumption of semi-algebraicity is not superfluous. In what follows, for a set $S$, the number of elements in $S$ will be denoted by $S^{\#}$. We begin with the following simple proposition.

\begin{prop}\label{prop:fin}
Let $F\colon \R^n\rightrightarrows\R^m$ be a semi-algebraic set-valued mapping whose graph has dimension no greater than $n$. Then there exists $\beta\in\mathbb{N}$ such that for a generic set of points $c\in \R^n$, we have $F(c)^{\#}\leq \beta$. 
\end{prop}
{\pf
Let $D = \mbox{\rm dom}\, F$. If the dimension of $D$ is strictly less than $n$, then we are done since then the complement of $D$ is a set satisfying the claimed property with $\beta=0$. Thus assume that $D$ has dimension $n$. Applying Corollary~\ref{cor:svh} to the mapping $F$, we get a finite partition of $D$ into semi-algebraic sets $\{C_i\}$, such that $$\mbox{\rm gph}\,F|_{C_i}\cong C_i\times F(c)$$ for any $c\in C_i$. Let $C_i$ be a partitioning set of maximal dimension. So $\dim C_i=n$ and we have $$n\geq\dim \mbox{\rm gph}\,F|_{C_i}=n+\dim F(c)$$ for any $c\in C_i$. Thus $\dim F(c)=0$ and since $F(c)$ is a semi-algebraic set, it must be finite. Since this argument holds for any $C_i$ of maximal dimension, we have that for a generic vector $c$, the set $F(c)$ is finite. Observe that if $F(c)$ is a finite non-empty set, then $F(c)^{\#}$ is equal to the number of connected components of $F(c)$. By Remark~\ref{cor:fiber}, there exists $\beta\in\mathbb{N}$ such that for all $c\in \R^n$, the number of connected components of $F(c)$ is no greater than $\beta$. So in particular, for generic $c$, we have $F(c)^{\#}\leq \beta$.}\qed   

\begin{cor}\label{cor:gen_unique}
Let $f\colon\R^n\rightarrow\overline{\R}$ be a semi-algebraic function and consider the collection of perturbed functions $h_v(x)=f(x)-\langle v,x\rangle$, parametrized by vectors $v\in\R^n$. Then there exists a positive integer $\beta$, such that for generic $v\in\R^n$, the number of Clarke-critical points of the perturbed function $h_v$ is no greater than $\beta$.
\end{cor}
{\pf
Observe $$0\in\partial_{c}h_v(x)\Leftrightarrow v\in\partial_{c}f(x)\Leftrightarrow x\in (\partial_{c}f)^{-1}(v).$$ Thus the set $(\partial_{c}f)^{-1}(v)$ is equal to the set of Clarke-critical points of the function $h_v$. By Theorem~\ref{thm:main}, we have $\dim \mbox{\rm gph}\, \partial_c f\leq n$, hence $\dim \mbox{\rm gph}\, (\partial_c f)^{-1}\leq n$. Applying Theorem~\ref{prop:fin} to $(\partial_c f)^{-1}$, we deduce that there exists a positive integer $\beta$, such that for generic $v$, we have $((\partial_c f)^{-1}(v))^{\#}\leq \beta$. The result follows.}\qed

\begin{cor}\label{cor:gen_nondeg}
Let $f\colon\R^n\rightarrow\overline{\R}$ be a semi-algebraic function and consider the collection of perturbed functions $h_v(x)=f(x)-\langle v, x\rangle$, parametrized by vectors $v\in\R^n$. Then for generic $v\in\R^n$, every Clarke-critical point of the function $h_v$ is nondegenerate. 
\end{cor}

Corollary~\ref{cor:gen_nondeg} follows immediately from the observation $$0 \in \mbox{ri}\, \partial_c h_v(x)\Leftrightarrow v\in \mbox{ri}\, \partial_c f(x),$$ and the following result. 
\begin{cor}\label{cor:Larman}
Let $f\colon\R^n\rightarrow\overline{\R}$ be a semi-algebraic function. Then for generic $v\in\R^n$, we have that $$x\in (\partial_c f)^{-1}(v) \Longrightarrow v\in \mbox{\rm ri}\, \partial_c f(x).$$
\end{cor}
{\pf 
Let $D = \mbox{\rm dom}\, \partial_c f$.
Consider the semi-algebraic set-valued mapping
$$\widetilde{F}\colon\R^n\rightrightarrows\R^n, ~~x\mapsto \textrm{rb }\partial_c f(x).$$

Our immediate goal is to show that the dimension of $\mbox{\rm gph}\,\widetilde{F}$ is no greater than $n-1$. Observe that for each $x\in\R^n$, we have $\widetilde{F}(x)\subset\partial_c f(x)$. Applying Corollary~\ref{cor:svh} to the mapping $\partial_c f$, we get a finite partition of $D$ into semi-algebraic sets $\{X_i\}$, such that $$\mbox{\rm gph}\,\partial_c f|_{X_i}\cong X_i\times \partial_c f(x)$$ and $$\mbox{\rm gph}\,\widetilde{F}|_{X_i}\cong X_i\times \widetilde{F}(x)$$ for any $x\in X_i$ (for each $i$). By Theorem~\ref{thm:main}, we have that $$n\geq\dim \mbox{\rm gph}\,\partial_c f|_{X_i} = \dim X_i+ \dim \partial_c f(x).$$    
Since $\widetilde{F}(x)=\textrm{rb }\partial_c f(x)$, it follows that $$\dim \widetilde{F}(x) \leq \dim \partial_c f(x)-1.$$ 
Therefore $$\dim \mbox{\rm gph}\,\widetilde{F}|_{X_i}= \dim X_i+ \dim \widetilde{F}(x)\leq \dim X_i+\dim \partial_c f(x)-1\leq n-1.$$
Thus $$\dim  \mbox{\rm gph}\, \widetilde{F}= \dim \Big(\bigcup_i \mbox{\rm gph}\, \widetilde{F}|_{X_i}\Big)\leq n-1.$$ And so if we let $$\pi\colon\mbox{\rm gph}\, \widetilde{F}\rightarrow\R^n$$ be the projection onto the last $n$ coordinates, we deduce that $\dim \pi(\mbox{\rm gph}\, \widetilde{F})\leq n-1$. Finally, observe 
\[
\pi(\mbox{gph}\,\widetilde F) ~=~ 
\Big\{ v \in \R^n : v \in \mbox{rb}\,\partial_c f(x)~ \mbox{for some}~ x \in \R^n \Big\},
\] and so the result follows.}\qed

\begin{rem}
{\rm Observe that if a convex function has finitely many minimizers then, in fact, it has a unique minimizer. Thus, for a proper convex semi-algebraic function, Corollaries~\ref{cor:gen_unique} and \ref{cor:gen_nondeg} reduce to Theorem~\ref{thm:gen_conv}.}  
\end{rem}

\begin{rem}\label{rem:lip}
{\rm In Corollaries~\ref{cor:gen_unique} and \ref{cor:gen_nondeg}, if the function $f$ is not semi-algebraic, then the results of these corollaries can fail. In fact, these results can fail even if the function $f$ is locally Lipschitz continuous. For instance, there is a locally Lipschitz function $f\colon\R\rightarrow\R$ such that $\partial_c f(x)=[-x,x]$ for every $x\in\R$. See the article of Borwein-Moors-Wang \cite{counter}. 
For all $v\in\R$, the perturbed function $h_v$ has infinitely many critical points, and for all $v\in\R\setminus\{0\}$, the function $h_v$ has critical points that are degenerate.}  
\end{rem}

\section{Composite Optimality Conditions} \label{sec:composite}
Consider a composite optimization problem $\displaystyle\min_x\, g(F(x))$. It is often computationally more convenient to replace the criticality condition $0\in\partial(g\circ F)(x)$ with the potentially different condition $0\in\nabla F(x)^{*}\partial g(F(x))$, related to the former condition by an appropriate chain rule. See for example the discussion of Lagrange multipliers in \cite{lag}. Thus it is interesting to study the graph of the set-valued mapping $x\mapsto \nabla F(x)^{*}\partial g(F(x))$.   
\subsection{Dimensional Analysis of the Chain Rule.}
The following is a standard result in subdifferential calculus.
\begin{thm}\cite[Theorem 10.6]{VA}\label{thm:chain}
Consider a function $g\colon\R^m\to\overline{\R}$ and a smooth mapping $F\colon\R^n\to\R^m$.
Then at any point $\bar{x}\in \R^n$, one has 
$$\hat{\partial} (g\circ F)(\bar{x}) \supset \nabla F(\bar{x})^{*}\hat{\partial}g(F(\bar{x})).$$
\end{thm}
Now assuming that the functions $g$ and $F$ in the theorem above are semi-algebraic, we immediately deduce, using Theorem~\ref{thm:main}, that the dimension of the graph of the mapping $x\mapsto\nabla F(\bar{x})^{*}\hat{\partial}g(F(\bar{x}))$ is at most $n$. 
One can ask what happens more generally in the case of the limiting and Clarke subdifferentials. It is well known that the inclusion $$\partial (g\circ F)(\bar{x}) \supset \nabla F(\bar{x})^{*}\partial g(F(\bar{x}))$$ is only guaranteed to hold under certain conditions \cite[Theorem 10.6]{VA}. The Clarke case is similar \cite[Theorem 2.3.10]{clarke}. Hence, a priori, the dimension of the graph of the set-valued mapping $x\mapsto \nabla F(x)^{*}\partial g(F(x))$ is unclear. In this section, we will show that if $g$ is lower semicontinuous, then this dimension is no greater than $n$ and we will derive some consequences.

The proofs of Proposition~\ref{prop:strat} and Theorem~\ref{thm:main} are self contained and purely geometric. There is, however, an alternative approach using \cite[Proposition 4]{Lewis-Clarke}, which will be useful for us. We state this proposition now. We denote the linear subspace of $\R^n$ parallel to a nonempty convex set $S\subset\R^n$ by $\para S$. 

\begin{prop}\cite[Proposition 4]{Lewis-Clarke}\label{prop:projection}
Consider a proper, lower semicontinuous, semi-algebraic function $g\colon\R^m\to\overline{\R}$. Then there exists a Whitney stratification $\{M_i\}$ of the domain of $g$ such that for each stratum $M_i$ and for any point $x\in M_i$, the inclusion $\para\partial_c g(x)\subset N_{M_i}(x)$ holds. 
\end{prop}

Before proceeding, we record the following special case of Theorem~\ref{thm:chain}. Consider a smooth function $F\colon\R^n\to\R^m$ and a nonempty set $Q\subset \R^m$. Consider any point $\bar{x}\in \R^n$. Applying Theorem~\ref{thm:chain} to the indicator function of $Q$, we deduce $$\hat{N}_{F^{-1}(Q)}(\bar{x})\supset \nabla F(\bar{x})^{*}\hat{N}_{Q}(F(\bar{x})).$$ If we let $Q=F(X)$, for some set $X\subset\R^n$, then we obtain 
\begin{equation}\label{eq:conv}
\hat{N}_{F^{-1}(F(X))}(\bar{x})\supset \nabla F(\bar{x})^{*}\hat{N}_{F(X)}(F(\bar{x})).
\end{equation} 

\begin{thm}\label{thm:calc}
Consider a proper, lower semicontinuous, semi-algebraic function $g\colon\R^m\to\overline{\R}$ and a smooth semi-algebraic mapping $F\colon\R^n\to\R^m$. Then the graph of the semi-algebraic set-valued mapping $x\mapsto \nabla F(x)^{*}\partial_c g(F(x))$  has dimension no greater than $n$. 
\end{thm}
{\pf 
Consider the Whitney stratification $\{M_i\}$ of $\dom g$ that is guaranteed to exist by applying Proposition~\ref{prop:projection} to the function $g$.
Now applying Theorem~\ref{theorem:strat} to the mapping $F$, we obtain a Whitney stratification $\{X_i\}$ of $\R^n$ and a Whitney stratification $\{K_j\}$ of $\R^m$ compatible with $\{M_i\}$ such that for each index $i$, we have $F(X_i)=K_j$ for some index $j$. Fix some stratum $X$ and a point $\bar{x}\in X$. If $F(X)$ is not a subset of the domain of $g$, then clearly $\nabla F(\cdot)^{*}\partial_c g(F(\cdot))|_X\equiv \emptyset$. Hence, we only consider $X$ such that $F(X)\subset\dom g$. Let $M$ be the stratum satisfying $F(X)\subset M$. Observe by our choice of the stratification $\{M_i\}$, we have $$\nabla F(\bar{x})^{*}\partial_c g(F(\bar{x}))\subset \nabla F(\bar{x})^{*}v+ \nabla F(\bar{x})^{*}N_M(F(\bar{x})),$$ for some vector $v\in\R^m$. Hence we have the inclusions
\begin{equation}\label{eq:main_ineq}
\para \nabla F(\bar{x})^{*}\partial_c g(F(\bar{x}))\subset \nabla F(\bar{x})^{*}N_M(F(\bar{x}))\subset \nabla F(\bar{x})^{*}N_{F(X)}(F(\bar{x})),
\end{equation}
where the last inclusion follows since the manifold $F(X)$ is a subset of $M$. Combining (\ref{eq:conv}) and (\ref{eq:main_ineq}), we obtain 
$$\para \nabla F(\bar{x})^{*}\partial_c g(F(\bar{x}))\subset \hat{N}_{F^{-1}(F(X))}(\bar{x})\subset N_X(\bar{x}),$$ where the last inclusion follows since the manifold $X$ is a subset of $F^{-1}(F(X))$. So we deduce $$\dim \nabla F(\bar{x})^{*}\partial_c g(F(\bar{x}))\leq n-\dim X.$$ Since the point $\bar{x}$ was arbitrarily chosen from $X$, we conclude, using Corollary~\ref{prop:const}, the inequality $\dim \gph\nabla F(\cdot)^{*}\partial_c g(F(\cdot))|_X\leq n$. Taking the union over the strata $\{X_i\}$ yields $$\dim \gph\nabla F(\cdot)^{*}\partial_c g(F(\cdot))\leq n,$$ as we claimed. 
}\qed

Observe that Theorem~\ref{thm:calc} is a generalization of Theorem~\ref{thm:main}. This can easily be seen by taking $F$ to be the identity map in Theorem~\ref{thm:calc}.

\begin{rem}
{\rm Consider a proper, lower semicontinuous, semi-algebraic function $g\colon\R^m\to\overline{\R}$ and a smooth semi-algebraic mapping $F\colon\R^n\to\R^m$ satisfying $\dom g\circ F= F^{-1}(\dom g)\neq\emptyset$. A natural question, in line with Theorem~\ref{thm:equation}, is whether the graph of the mapping $x\mapsto \nabla F(x)^{*}\partial_c g(F(x))$ has dimension exactly $n$. In fact, there is no hope for that to hold generally. For instance, it is possible to have $\partial_c g(y)=\emptyset$ for every point $y$ in the image of $F$. This example, however motivates the following easy proposition.

\begin{prop}
Consider a proper, lower semicontinuous, semi-algebraic function $g\colon\R^m\to\overline{\R}$ and a smooth semi-algebraic mapping $F\colon\R^n\to\R^m$.
Assume that the set $F^{-1}(\dom \hat{\partial}g)$ has a nonempty interior. Then the graph of the set-valued mapping $\nabla F(\cdot)^{*}\hat{\partial} g(F(\cdot))$ has dimension exactly $n$. Analogous results hold in the limiting and Clarke cases.
\end{prop}
{\pf 
We show the theorem only for the case of the regular subdifferential. The proof remains unchanged for the other cases.
Clearly, as a result of Theorem~\ref{thm:calc}, it is sufficient to show that the inequality, $\dim \gph \nabla F(\cdot)^{*}\hat{\partial} g(F(\cdot))\geq n$, holds. To that effect, consider an open set $N$ contained in the interior of the set $F^{-1}(\dom \hat{\partial}g)$. Then for any point $x\in N$, the set $\hat{\partial} g(F(x))$ is nonempty. Hence, the the map $x\mapsto \nabla F(x)^{*}\hat{\partial} g(F(x))$ has nonempty values on $N$. In particular, letting $\pi\colon\R^n\times\R^n\to\R^n$ be the projection onto the first $n$ coordinates, we see $$\pi(\gph \nabla F(\cdot)^{*}\hat{\partial} g(F(\cdot))|_N)=N.$$ Hence we conclude $\dim \gph \nabla F(\cdot)^{*}\hat{\partial} g(F(\cdot))|_N\geq \dim N=n$, thus completing the proof.     
}\qed

The following easy result, which we state without proof, gives a different approach.

\begin{prop}
Consider a proper, continuous, semi-algebraic function $g\colon\R^m\to\overline{\R}$ and a smooth semi-algebraic mapping $F\colon\R^n\to\R^m$. If the constraint qualification of \cite[Theorem 10.6]{VA} holds at some point $x\in F^{-1}(\dom g)$, then the graph of the set-valued mapping $\nabla F(\cdot)^{*}\partial g(F(\cdot))$ has dimension exactly $n$. An analogous statement holds in the Clarke case.   
\end{prop}
}
\end{rem}


\subsection{Consequences}\label{sec:cons}
Let $F\colon\R^n\to\R^m$ be a smooth mapping and $g\colon\R^m\to\overline{\R}$ a proper lower semicontinuous function. (For simplicity, here we assume that the mapping $F$ is defined on all of $\R^n$. However the whole section extends immediately to a mapping $F$ defined only on an open subset $U\subset\R^n$.) Consider the following collection of composite minimization problems, parametrized by vectors $v\in\R^n$. $$(P(v))\qquad\min_{x\in\R^n} g(F(x))-\langle v,x \rangle$$ For a point $\bar{x}$ to be a minimizer for $P(v)$, the inclusion $v\in\partial(g\circ F)(\bar{x})$ must necessarily hold. As discussed in the beginning of the section, it is often more convenient to replace this condition with the potentially different condition $v\in\nabla F(\bar{x})^{*}\partial g(F(\bar{x}))$. This motivates the following definition.
\begin{defn}
{\rm
We say that a point $x$ is {\em Clarke critical} for the problem $(P(v))$ if the inclusion $v\in\nabla F(x)^{*}\partial_c g(F(x))$ holds, and we call such a critical point $x$  {\em nondegenerate} for the problem $(P(v))$ if the stronger property $v\in\ri \nabla F(x)^{*}\partial_c g(F(x))$ holds.
}
\end{defn}

We are now in position to state a natural generalization of Corollaries~\ref{cor:gen_unique} and \ref{cor:gen_nondeg}. 

\begin{cor}\label{cor:comp_unique2}
Let $F\colon\R^n\to\R^m$ be a semi-algebraic smooth function and $g\colon\R^m\to\overline{\R}$ a proper lower semicontinuous semi-algebraic function. Consider the following collection of optimization problems, parametrized by vectors $v\in\R^n$.
$$(P(v))\qquad \min_{x\in\R^n} g(F(x))-\langle v, x\rangle$$
Then there exists a positive integer $\beta$, such that for a generic vector $v\in\R^n$, the number of Clarke-critical points for the problem $(P(v))$ is no greater than $\beta$. Furthermore, for a generic vector $v\in\R^n$, every Clarke-critical point for the problem $(P(v))$ is nondegenerate. 
\end{cor}
{\pf Observe that by Theorem~\ref{thm:calc}, the graph of the mapping $x\mapsto \nabla F(x)^{*}\partial_c g(F(x))$  has dimension no greater than $n$. The proof now proceeds along the same lines as the proofs of Corollaries~\ref{cor:gen_unique} and \ref{cor:gen_nondeg}.
}\qed

Observe that Corollaries~\ref{cor:gen_unique} and \ref{cor:gen_nondeg} can be considered as special cases of Corollary~\ref{cor:comp_unique2}, in which the map $F$ is the identity map.

A noteworthy illustration of Corollary~\ref{cor:comp_unique2} is the problem of constrained minimization, which we discuss now. Let $f\colon\R^n\to\overline{\R}$ be a semi-algebraic function and $D\subset\R^n$ a closed semi-algebraic set. Consider the following collection of constrained minimization problems, parametrized by vectors $v\in\R^n$.
\begin{eqnarray*}
(P'(v))\qquad &\min &f(x)-\langle v,x \rangle\\
&\st  &x\in D
\end{eqnarray*}
Observe that $(P'(v))$ is equivalent to the problem $\displaystyle\min_{x\in\R^n} g(F(x))-\langle v, x\rangle$, where we define $F(x)=(x,x)$ and $g(x,y)=f(x)+\delta_D(y)$. 

Hence, in the sense of composite minimization, it is easy to check that a point $x\in D$ is {\em Clarke critical} for the problem $(P'(v))$ if $v\in\partial_c f(x)+N_D^c(x)$, and  such a critical point $x$  is {\em nondegenerate} for the problem $(P'(v))$ if the stronger property $v\in\ri \partial_c f(x)+\ri N_D^c(x)$ holds. 



\begin{cor}\label{cor:gen_unique2}
Let $f\colon\R^n\to\R$ be a semi-algebraic function and let $D$ be a closed, semi-algebraic set. Consider the following collection of optimization problems, parametrized by vectors $v\in\R^n$.
\begin{eqnarray*}
(P(v))\qquad &\min &f(x)-\langle v, x\rangle\\
&\st  &x\in D
\end{eqnarray*}
Then there exists a positive integer $\beta$, such that for a generic vector $v\in\R^n$, the number of Clarke-critical points for the problem $(P(v))$ is no greater than $\beta$. Furthermore, for a generic vector $v\in\R^n$, every Clarke-critical point for the problem $(P(v))$ is nondegenerate. 
\end{cor}
{\pf This follows directly from Corollary~\ref{cor:comp_unique2}.
}\qed

\begin{acknowledgements}
Thanks to Alex D. Ioffe for suggesting the extension we pursue in Section~\ref{sec:composite}.
\end{acknowledgements}

\bibliographystyle{hplain}
\small
\parsep 0pt
\bibliography{dim_graph}

\end{document}